\documentclass{amsart}
\usepackage{graphicx,euscript,eufrak,amsmath,amsthm, color, mathrsfs,hyperref}
\usepackage[margin=1.3in]{geometry} 

\newtheorem{theorem}{Theorem}[section]
\newtheorem{lemma}{Lemma}[section]
\newtheorem{corollary}[theorem]{Corollary}
\theoremstyle{definition}
\newtheorem{remark}{Remark}[section]

\newcommand{\ZZ}{\mathbb{Z}}
\newcommand{\kk}{\kappa}
\newcommand{\vk}{\vec{k}}

\newcommand{\ve}{\vec{e}}
\newcommand{\al}{\alpha}

\newcommand{\vE}{\vec{E}}
\newcommand{\vD}{\vec{D}}

\newcommand{\vn}{\vec{n}}
\newcommand{\dd}{\delta}

\newcommand{\ti}{\textit}

\long\def\ignore#1{}

\def\nd{\noindent}

\def\dd{\delta}

\def\ti{\textit}
\def\la{\left\langle}
\def\ra{\right\rangle} 
\def\lb{\left\lbrace}
\def\rb{\right\rbrace}

\def\RR{\mathbb{R}}
\def\ZZ{\mathbb{Z}}

\def\vx{\vec{x}}
\def\vy{\vec{y}}

\begin{document}

\title{The Residual Set Dimension of a Generalized Apollonian Packing}

\author{Daniel Lautzenheiser}
\begin{abstract}
	We view space-filling circle packings as subsets of the boundary of hyperbolic space subject to symmetry conditions based on a discrete group of isometries. This allows for the application of counting methods which admit rigorous upper and lower bounds on the Hausdorff dimension of the residual set of a generalized Apollonian circle packing. This dimension (which also coincides with a critical exponent) is strictly greater than that of the Apollonian packing.
\end{abstract}
\address{Eastern Sierra College Center, 4090 W. Line St, Bishop, CA, 93514}
\thanks{\nd Date: \today.}

\email{daniel.lautzenheiser@cerrocoso.edu}

\maketitle

\section{introduction}

The Apollonian packing is an infinite collection of circles whose corresponding discs cover a region almost everywhere. The residual set (the space not covered) is a Cantor-like set with Hausdorff dimension satisfying
\begin{align}\label{apollonian new bound}
    1.302327 < \delta_\mathcal{A} < 1.310876.
\end{align}
The main result of this paper is that the Hausdorff dimension of a non-Apollonian packing to be introduced briefly satisfies
\begin{align}
    1.327266 < \delta < 1.348771.
\end{align}
This is the first set of rigorous bounds on the residual set dimension for a packing other than the Apollonian packing. 

Given a connected open region $\Omega \subseteq \RR^2$, we will refer to a \ti{circle packing}, or just \ti{packing}, $\mathcal{P}$, as a disjoint union of open discs in $\Omega$. The \ti{residual set} is $\mathcal{R} = \Omega - \mathcal{P}$. We say that a given packing is \ti{complete} if $\mathrm{Vol}(\mathcal{R}) = 0$. We will call a packing \ti{bounded} if $\Omega$ is a bounded subset of $\RR^2$. Associated to a bounded circle packing is a sequence of radii $ \{r_n\}_{n=1}^\infty $ which can be reordered to decrease to zero. Consider the ``L-function"     
\begin{equation}\label{formal radius series}
L(t) = \sum_{n=1}^\infty r_n^t.
\end{equation}
Notice that $ L(2) = \frac{1}{\pi} \mathrm{Vol}(\mathcal{P}) $. If $ t = 1 $, the sequence of partial sums represents a sum of circumferences apart from a factor of $ \frac{1}{2 \pi} $. For the Apollonian packing, this full series diverges \cite{Mergelyan1952uniform,wesler1960infinite}. It is then natural to study the \ti{critical exponent}
\begin{align}\label{critical exponent}
\begin{split}
S =  \inf \lb t : \sum_{n=1}^\infty r_n^t < \infty  \rb = \sup \lb t : \sum_{n=1}^\infty r_n^t = \infty  \rb.
\end{split}
\end{align}
For the Apollonian and generalized (or non-Apollonian) packing to be described briefly, the Hausdorff dimension of the residual set coincides with the critical exponent \cite{kontorovich2011apollonian}.

It has been known for some time that there are non-Euclidean aspects to packing problems \cite{Maxwell1981Space}. David Boyd's ``separation" formula \cite{Boyd1973Separation} (which can be traced back to Darboux \cite{Darboux1872Points} and Clifford \cite{Clifford1882Powers}) is used to determine a \ti{polyspherical} coordinate system for packed spheres. With the machinery of hyperbolic (or \ti{Lobackevsky}) geometry, we may interpret the separation of two spheres as the minimal distance connecting two disjoint hyperbolic planes. We will use the \ti{pseudosphere} or \ti{vector model} of hyperbolic geometry. See \cite{Baragar2001Survey,Ratcliffe2006Foundations,Dolgachev2016Orbital,Apanasov2011Conformal} for additional introductory information.

Lorentz space $ \mathbb{R}^{n-1,1} $ is defined as the set of vectors in $ \mathbb{R}^n $ together with the \ti{Lorentz product}
\begin{align}\label{Lorentz product (unweighted)}
\vx \circ \vy = x_1y_1 + \cdots + x_{n-1}y_{n-1} - x_ny_n.
\end{align}  
More generally, any symmetric $ n \times n $ matrix $ J $ of signature $ (n-1,1) $ will define, up to a change of basis, a Lorentz product via $ \vx \circ \vy = \vx^tJ\vy $. The surface cut by the equation $ \vx \circ \vx = -1 $ is a hyperboloid of two sheets. The top or forward sheet
\[ \mathcal{H} = \lb \vx : \quad \vx \circ \vx = -1, \quad x_n>0  \rb  \]
together with the metric $ d $ where
\begin{equation}\label{distance_metric}
\vx \circ \vy = ||\vx|| ||\vy|| \cosh d(\vx,\vy) = - \cosh d(\vx,\vy) 
\end{equation}
is a model of hyperbolic geometry, $ \mathbb{H}^{n-1} $. The linear maps preserving the Lorentz product, called $ \ti{Lorentz transformations} $ can be identified with the infinite matrix group
\begin{align}
\begin{split}
\mathcal{O}_J(\RR) = \{ T\in{M_{n \times n}(\RR)}: T \vx \circ T \vy =  \vx \circ \vy \quad \text{ for all }  \vx, \vy \in{\RR^n}  \}, 
\end{split}
\end{align}
and the group of isometries of the model $ \mathcal{H} $ can be identified with 
\begin{equation}
\mathcal{O}_J^+(\RR) = \{ T\in{\mathcal{O}_J(\RR)}: T \mathcal{H} = \mathcal{H}  \}.
\end{equation}
Reflection through $\vec{n} \circ \vx = 0$ is
\begin{align}
    R_{\vec{n}}(\vec{x}) = \vec{x} - 2 \cdot \mathrm{Proj}_{\vec{n}}(\vx) = \vx - 2 \frac{\vx \circ \vn}{\vn \circ \vn}\vn. 
\end{align}
A subgroup of interest is the discrete group $ \mathcal{O}_J^+(\ZZ) $ and, generally speaking, the symmetries of a packing (often reflections) may be represented by a subgroup of finite index in $ \mathcal{O}_J^+(\ZZ) $.
Let $ n=4 $ and fix the standard basis $ \{\ve_1,\ve_2,\ve_3,\ve_4 \} $. The matrix
\begin{equation}
J = (\ve_i \circ \ve_j) = -\begin{pmatrix}
-1 & 1 & 1 & 1 \\
1 & -1 & 1 & 1 \\
1 & 1 & -1 & 1 \\
1 & 1 & 1 & -1
\end{pmatrix}
\end{equation}
diagonalizes with eigenvalues $ \{2,2,2,-2\} $. Thus, the bilinear form
$(\vx,\vy) \mapsto \vx^t J \vy = \vx \circ \vy $
provides a model of Lorentzian $ 4- $space, which is isometrically equivalent to the Poincar{\'e} upper half space model \cite{Baragar2017Higher}. With the interpretation that $ \ve_j $ represents a normal vector to the plane $ \vx \circ \ve_j = 0 $, then in view of formula (\ref{distance_metric}), $ J $ conveys to us that the configuration of basis vectors represents mutually tangent planes. Since planes intersect the boundary $ \partial \mathbb{H}^3 $ in Euclidean circles or planes, $ J $ defines a configuration of 4 mutually tangent circles on the boundary. It is then, in general, a non-trivial problem to identify a subgroup $ \Gamma \leq \mathcal{O}_J^+(\ZZ) $ which acts on these four planes (or faces) whose orbit is a complete packing (in this case, the Apollonian packing.) 

\begin{figure}[h]
	\begin{center}
		\includegraphics[scale=0.15]{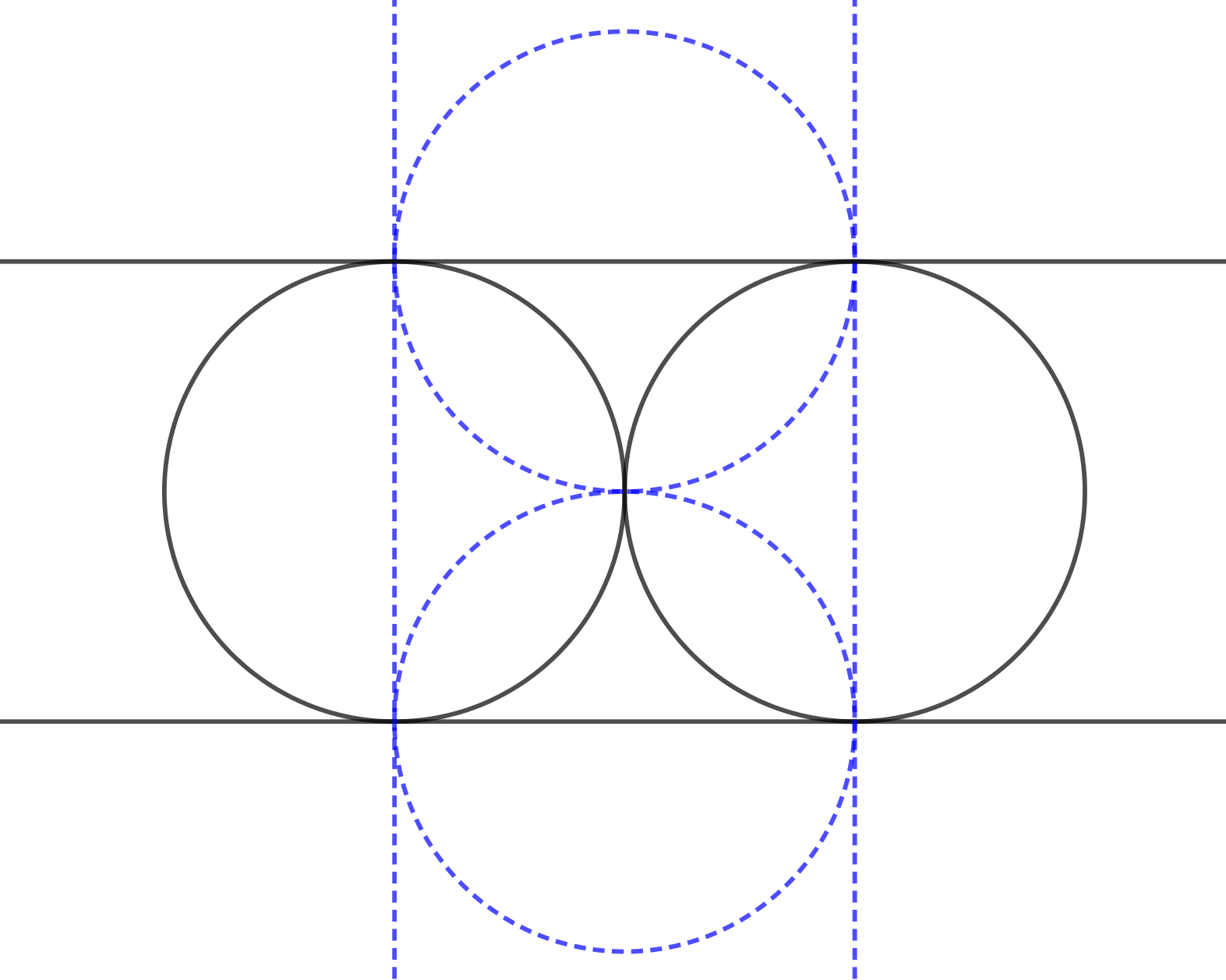}
		\includegraphics[scale=0.35]{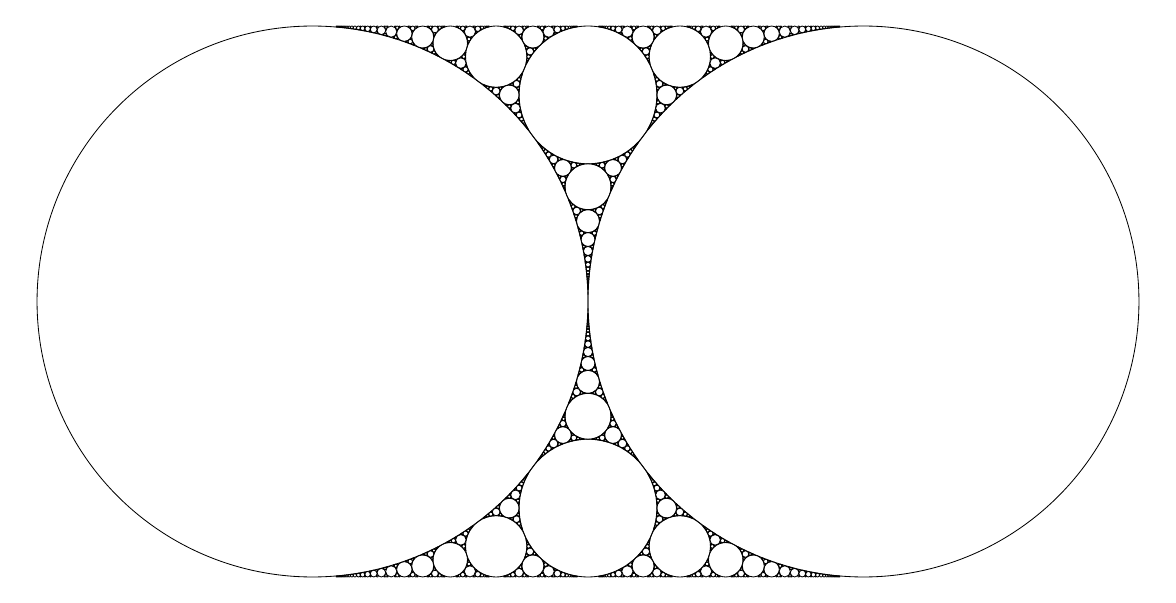}
	\end{center}
	\caption{An initial configuration of the four basis vectors along with some symmetries (dotted lines) as visualized on the boundary of the model of hyperbolic space induced by the circulant $ J $ above. $ \ve_1 $ and $ \ve_2 $ represent the two large circles while $ \ve_3 $ and $ \ve_4 $ the horizontal lines. \label{apollonian packing initial config pic}}
\end{figure}

\section{The Separation-3 (or Boyd/Mallows) Packing}
In the case of the Apollonian packing, the basis vectors $ \{\ve_1, \ve_2,\ve_3,\ve_4  \} $ may be chosen to represent mutually tangent circles, which in turn implies that the separation matrix has the above circulant form. Boyd \cite{Boyd1974New} showed the existence of several packings which do not possess the mutually tangent condition. One such set of examples begins if the disks represented by $ \ve_1 $ and $ \ve_2 $ are disjoint \cite{baragar2020game}. Let $r_1$ and $r_2$ be the radii of $\vec{e}_1$ and $\vec{e}_2$ on the boundary and $s$ the distance between their centers. The separation formula \cite{Boyd1973Separation} 
\begin{equation}\label{separation formula}
\cosh d(\vec{e}_1,\vec{e}_2) = \frac{s^2-r_1^2-r_2^2}{2r_1r_2}
\end{equation}
implies that if $ r_1 = r_2 = \frac{1}{2} $ and $ s = \sqrt{2} $, then we get a separation of $ 3 $; and via (\ref{distance_metric}), $- (\vec{e_1} \circ \vec{e_2}) = 3$. This leads to using the matrix
\[ J = -\begin{pmatrix}
-1&3&1&1\\3&-1&1&1\\1&1&-1&1\\1&1&1&-1\\
\end{pmatrix} \]
to define a Lorentz product as well as a packing of disks with some distinct differences compared to the Apollonian packing. It is also possible to choose $ -J $ above and all the results are the same modulo minus signs. The existence of this particular packing was first shown in \cite{Boyd1974New}. Although Boyd provides the separation matrix and a general description of the construction process, it was likely not until \cite{Manna1991Precise} that pictures of this packing were published. Additional analysis was done in \cite{Mallows2010Generalization}. 

\begin{figure}[h]
	\begin{center}
		\includegraphics[scale = 0.3]{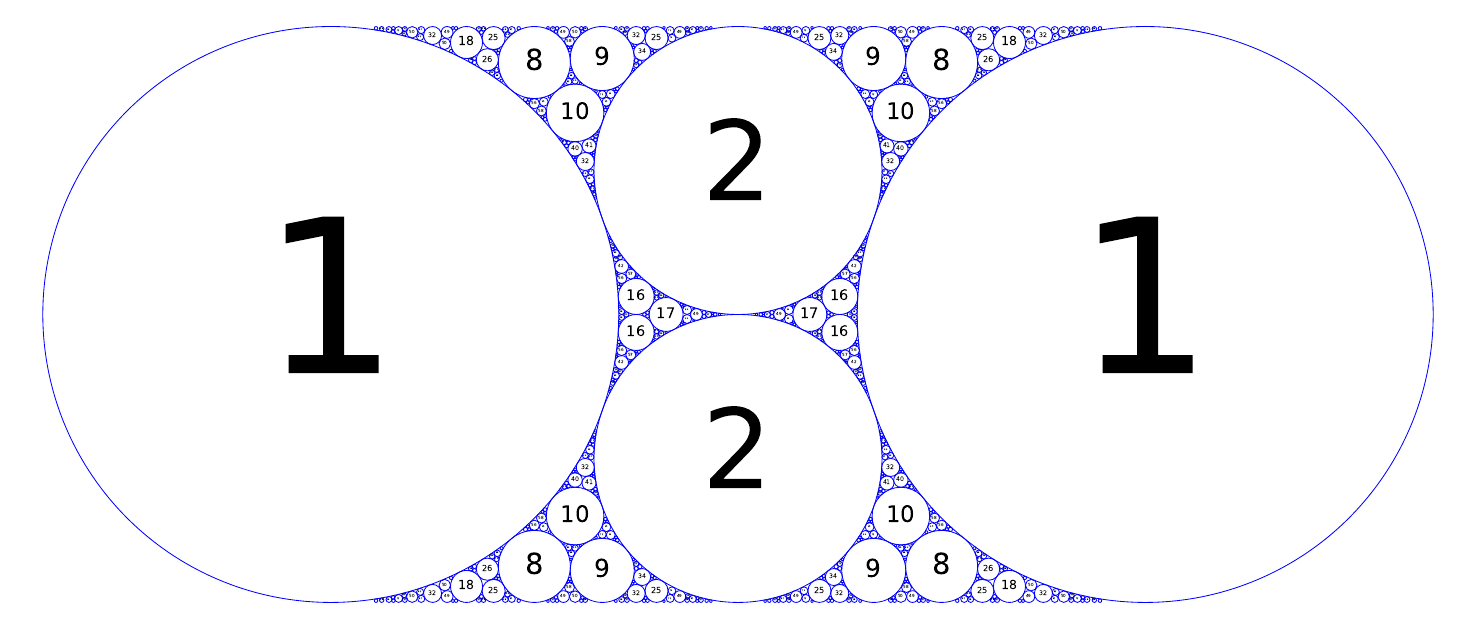}
	\end{center}
	\caption{A region of a Boyd/Mallows packing with integer curvatures.}
	\label{boyd mallows packing bounded}
\end{figure}

\begin{figure}[h]
	\begin{center}
		\includegraphics[scale = 0.15]{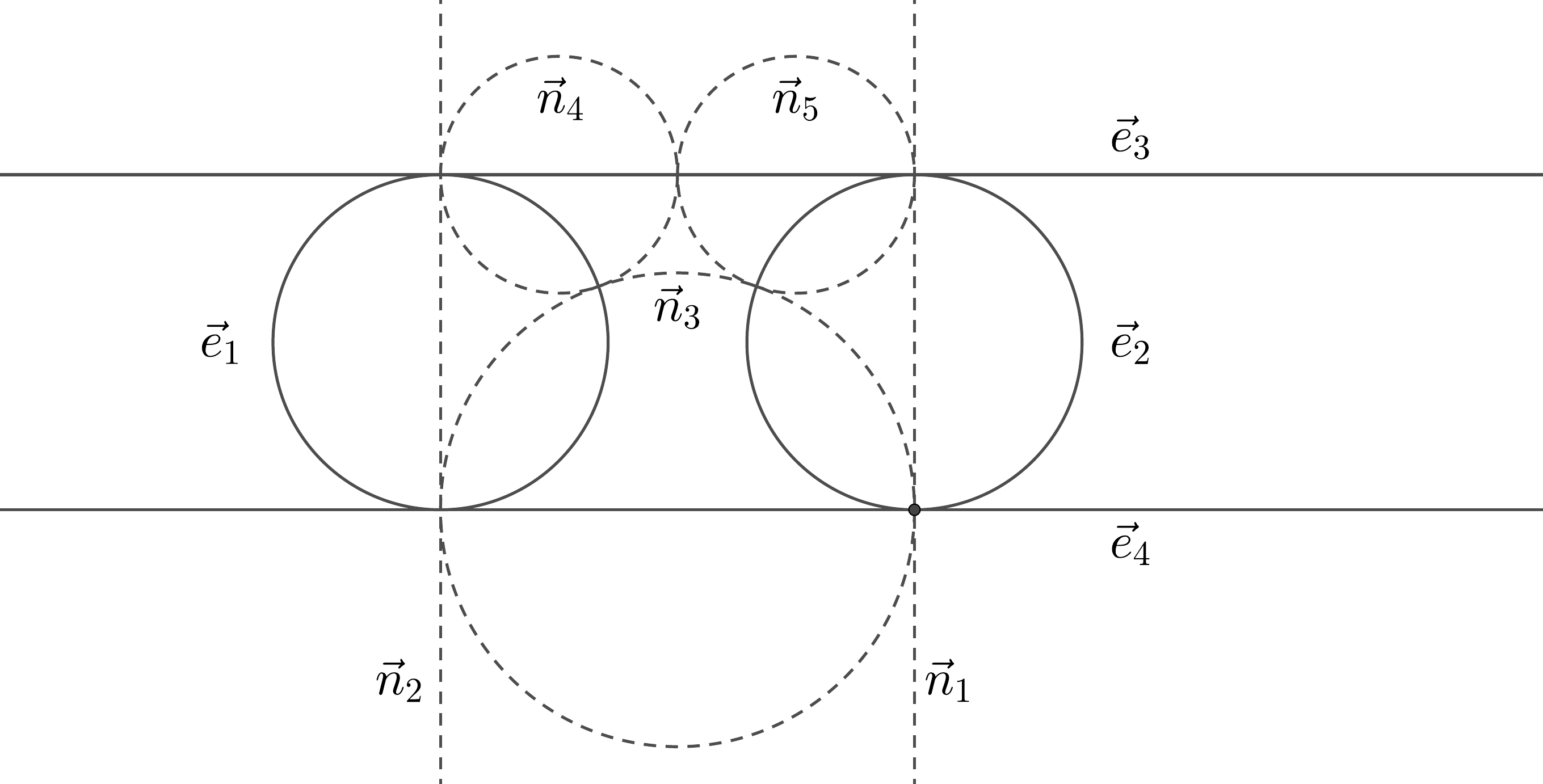}
	\end{center}
	\caption{A set of generating symmetries (dotted lines/circles) for the Boyd/Mallows packing which intersect tangentially.}
	\label{b/m symmetries}
\end{figure}

This packing may be generated by taking the image of   
\begin{align*}
  \Gamma = \la R_{\vn_1},R_{\vn_2},R_{\vn_3},R_{\vn_4},R_{\vn_5}  \ra  
\end{align*}
acting on the faces 
\begin{align*}
 \{ \ve_1,\ve_2,\ve_3,\ve_4,\ve_1+\ve_2-\ve_4 \}   
\end{align*}
where
$\vec{n}_1 = (-1,1,2,2)$, $\vec{n_2} = (1,-1,2,2)$, $\vec{n}_3 = (1,1,-2,0)$, $\vec{n}_4 = (3,1,2,-2)$, and $\vec{n}_5 = (1,3,2,-2)$.

\begin{figure}[h]
	\begin{center}
		\includegraphics[scale = 0.4]{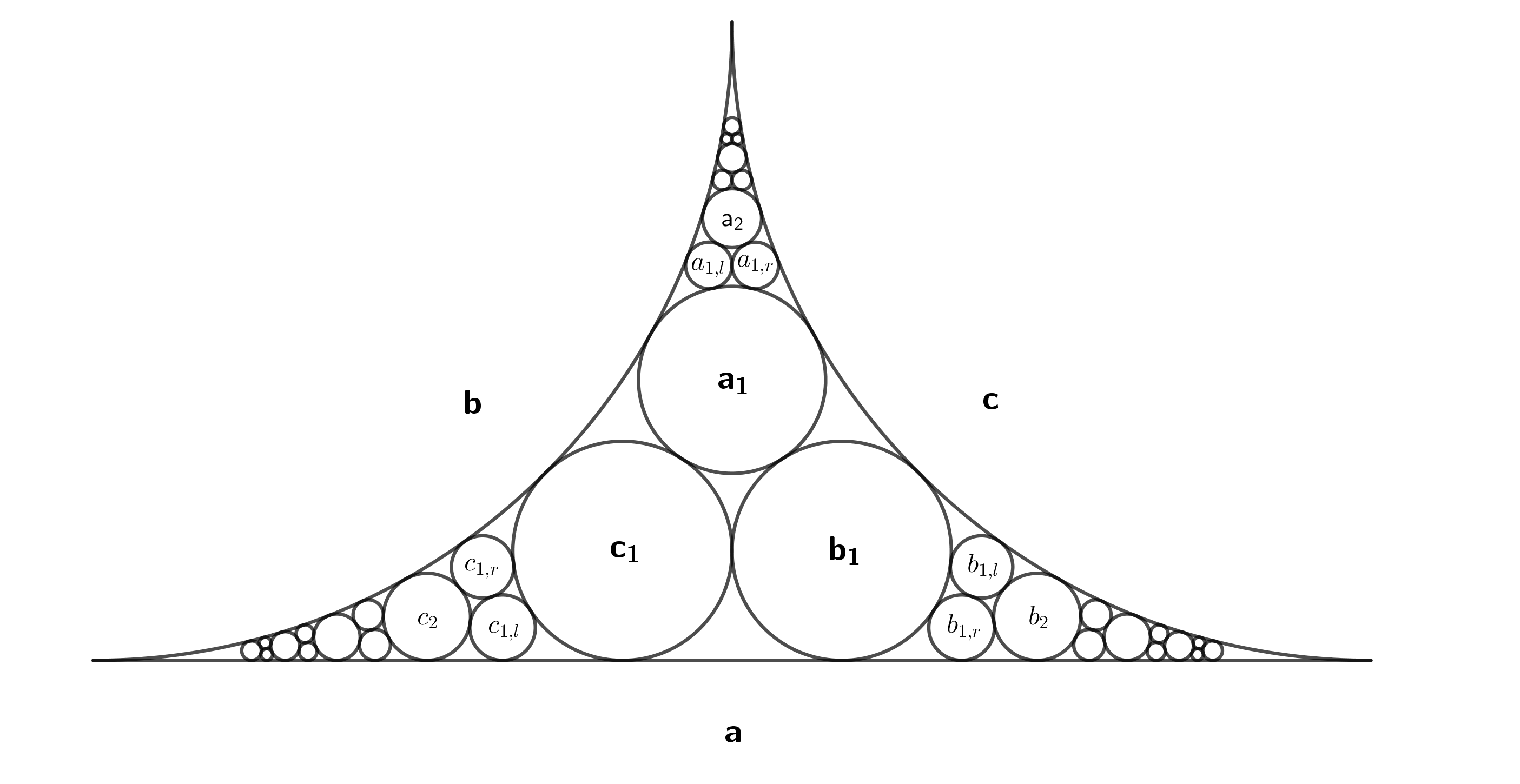}
	\end{center}
	\caption{The $ 0^{\text{th}} $ iteration necklace of the Boyd/Mallows packing.}
	\label{mallows necklace}
\end{figure}

Let $ \{r_n\} $ be a sequence of radii of the disks packed within a bounded region of a Boyd/Mallows or separation-3 packing. We wish to bound the new critical exponent $ S_{\text{Boyd/Mallows}} = \delta_{\text{Boyd/Mallows}} = \delta $ defined by (\ref{critical exponent}), which coincides with the Hausdorff dimension of the residual set \cite{Sullivan1984Entropy}. Set $ \vk = (a,b,c,d) = (\ve_1 \circ \vE, \ve_2 \circ \vE, \ve_3 \circ \vE, \ve_4 \circ \vE)  $, with $\vec{E}$ the point at infinity. Writing $ \vk^t J^{-1} \vk = 0 $ gives an analogue of the Descarte's quadruple relationship
\begin{equation}\label{Boyd/Mallows quadruple relationship}
-\frac{1}{8}(a^2+b^2)-\frac{1}{2}(c^2+d^2)+\frac{1}{4}ab+\frac{1}{2}(a+b)(c+d)= 0.
\end{equation}
Define the quadratic form
\[ K(\vk) = \vk^t J^{-1}\vk \]
and define $ c_n $ as being the smaller of the two solutions to $ K(c_n,c_{n-1},a,b) = 0 $. The sequence $ \{c_n\} $ is the set of curvatures of ``center disks" down the necklace opposite the triangle from the bounding circle $ C $ of curvature $ c=c_0 $ (see Figure \ref{mallows necklace}). By symmetry, $ \{ b_n \} $ and $ \{a_n\} $ are the curvatures of the center tails of circles in necklaces opposite bounding circles $ B $ and $ A $, respectively. To solve for $ a_1 $, we know that the the disks corresponding to $ a $ and $ a_1 $ are disjoint, so solving $ K(a_1,a,b,c) = 0 $ results in $ a_1 = a+2b+2c\pm \sqrt{8}\sqrt{ab+ac+bc} $.
Similarly, $ b_1 = 2a+b+2c \pm \sqrt{8}\sqrt{ab+ac+bc} $, and also, $ c_1 = 2a+2b+c\pm \sqrt{8}\sqrt{ab+ac+bc} $.

\begin{lemma} \label{melzak type lemma}
	Let $ A $, $ B $, and $ C $ be three pairwise externally tangent circles with curvatures $ a,b, $ and $ c $. Let $ \{C_n \} $ be the sequence of disks in which $ C_1 $ is the smaller of the disks tangent to $ A $ and $ B $, and at a separation of $ 3 $ from $ C $. Let $ C_n $ be the smaller of the disks tangent to $ A $ and $ B $ at a separation of $ 3 $ from $ C_{n-1} $. Then, for all $ n\in{\mathbb{N}} $, 
	\begin{align}\label{c_n or g_n function}
	c_n = c+2n^2(a+b)+n\sqrt{8}d
	\end{align}
	where $ d = \sqrt{ab+ac+bc} $. 
\end{lemma}

\begin{proof}:
	Let $ \{ \ve_1,\ve_2,\ve_3,\ve_4 \} $ be the standard basis with the above separation matrix and let 
	\begin{align*}
	\vk & = (\text{curv}(\ve_1), \text{curv}(\ve_2),\text{curv}(\ve_3),\text{curv}(\ve_4) \\
	&= (\ve_1 \circ \vE, \ve_2 \circ \vE, \ve_3 \circ \vE, \ve_4 \circ \vE) = (c,c_1,a,b)
	\end{align*}
	Consider the reflections $ R_1 = R_{(1,-1,2,2)} $ and $ R_2 = R_{(1,-1,0,0)} $. Then, $ R_2R_1 = P $ is a parabolic translation moving $ \ve_1 $ to $ \ve_2 $ and
	\[ P^n = \begin{pmatrix}
	1-n&-n&0&0\\n&n+1&0&0\\2n(n-1)&2n(n+1)&1&0\\2n(n-1)&2n(n+1)&0&1\\
	\end{pmatrix}. \]
	Then,  
	\begin{align*}
	c_n & =  \text{curv}(P^n\ve_1) = \text{curv}(1-n,n,2n(n-1),2n(n-1)) \\
	& = (1-n) (\ve_1 \circ \vE)  + n(\ve_2 \circ \vE) + 2n(n-1)((\ve_3 + \ve_4) \circ \vE) \\
	& = (1-n)c + nc_1+2n(n-1)(a+b).
	\end{align*}
	Note that, based on the picture determined by the curvatures $ a,b,c,c_1 $, then $ \vE$ (the point at infinity) is not $ \ve_3 + \ve_4 $.
	Using $ c_1 = 2a+2b+c+\sqrt{8}d $ gives the desired result.
\end{proof}
For example, using $ T(0,1,2) $ pictured above, we have $ c_1 = 8 $, $ c_2 = 18 $, $ c_3 = 32 $ etc. Since the necklace opposite the curvilinear triangle from either $ A $, $ B $ or $ C $ has three ``tails", we need formulas for the curvatures of the circles in the ``left" and ``right" tails.

\begin{lemma}\label{left necklace}
	Let $ A $, $ B $, and $ C $ be three pairwise externally tangent circles with curvatures $ a,b, $ and $ c $. Let $ \{C_{n,r} \} $ be the sequence of circles in which $ C_{1,r} $ is the smaller of the circles tangent to $ C_1 $ and $ C_2 $, and at a separation of $ 3 $ from $ A $. Let $ C_{n,l} $ be the smaller of the circles tangent to $ C_{n+1} $ and $ C_n $  at a separation of $ 3 $ from $ B $. Then, for all $ n\in{\mathbb{N}} $, 
	\begin{align}\label{c_{n,r} function}
	c_{n,r} = a + 2b+2c+4(n^2+n)(a+b)+(2n+1)\sqrt{8}d
	\end{align}
	and
	\begin{align}\label{c_{n,l} formula}
	c_{n,l} = 2a+b+2c+4(n^2+n)(a+b)+(2n+1)\sqrt{8}d
	\end{align}
	where $ d = \sqrt{ab+ac+bc} $.
\end{lemma}
\begin{proof}:
	As before, let $ \{ \ve_1,\ve_2,\ve_3,\ve_4 \} $ be a basis with the separation matrix $ J $ above and
	$\vk = (\ve_1 \circ \vE, \ve_2 \circ \vE, \ve_3 \circ \vE, \ve_4 \circ \vE) = (c,c_1,a,b) $. Consider the reflection $ R_3 = R_{(1,1,-2,0)} $. Then, $ R_3(\ve_3) = (1,1,-1,0) $. Using again the parabolic translation $P$ above, 
	\[ P^n(1,1,-1,0) = (1-2n,1+2n,4n^2-1,4n^2) \]
	implying
    \[ c_{n,r} = (1-2n)c + (1+2n)c_1+(4n^2-1)a+4n^2b. \]
	Substituting in $ c_1 = 2a+2b+c+\sqrt{8}\sqrt{ab+ac+bc} $ gives the desired result for $c_{n,r}$. Interchanging $ a $ and $ b $ gives $ c_{n,l} $. 
\end{proof}

Let $ T(a,b,c) $ be the curvilinear (or triply asymptotic) triangle bounded by three mutually externally tangent circles $ A,B,C $ of curvatures $ a,b,c $, with $ 0 \leq a \leq b \leq c $, and $ b>0 $. The condition $ b>0 $ guarantees that $ T(a,b,c) $ has finite area even if $ a = 0 $, in which case $ A $ is a line. For $ t>0 $, define 

\begin{equation} \label{M function B/M packing}
M(a,b,c;t) = \sum_{n=1}^\infty r_n^t
\end{equation}
where the $ r_n $ are the radii of the disks in the Boyd/Mallows packing within $ T(a,b,c) $ and the equality holds in the extended sense. To save writing space, we may suppress the variable $ t $, writing $ M(a,b,c) $. First note that based on the symmetries $\Gamma$, $ M $ is symmetric in the three variables $ a,b,c $.

\begin{lemma} \label{decreasing lemma B/M}
	$ M(a,b,c;t) $ is decreasing in each variable $ a,b,c $. $ M $ is strictly decreasing if  $ t > \dd $, and $ M(a,b,c;t)= + \infty $ if $ t < \dd$, where $\delta$ is the critical exponent.
\end{lemma}
Note that, apriori, $ \dd $ may depend on $ (a,b,c) $. We will see shortly that this is not the case.
\begin{proof}: If $ t < \dd $, equality holds trivially in the extended sense $ (\infty = \infty) $. Next let $ t > \dd$. Interchanging $a$ with $c$ or $b$ with $c$ gives, by symmetry, similar expressions for necklace curvatures $a_n,a_{n,l},a_{n,r}$, $b_n,b_{n,l}$, and $ b_{n,r}$. Let $ \epsilon > 0  $. If we replace $ c \mapsto c+ \epsilon $ then by Lemmas \ref{melzak type lemma} and \ref{left necklace}, $ a_n,a_{n,l},a_{n,r},b_n,b_{n,l},b_{n,r},c_n,c_{n,l},c_{n,r}$ strictly increase. We now expand $ M(a,b,c) $ by writing $ M(a,b,c) $ as a sum of 3 necklaces (each with 3 bands of circles), four central triangles, and 18 sub-triangle bands (6 from each necklace): 
\begin{align} \label{functional type equation 0th iteration sum b/m}
M(a,b,c) & = \sum_{n=1}^\infty \left( a_n^{-t} + a_{n,l}^{-t} + a_{n,r}^{-t} + b_n^{-t} + b_{n,l}^{-t} + b_{n,r}^{-t} + c_n^{-t} + c_{n,l}^{-t} + c_{n,r}^{-t}  \right) \nonumber \\
& + M(a,c_1,b_1)+M(b,c_1,a_1)+M(c,b_1,a_1)+M(c_1,b_1,a_1) \nonumber \\
& + \sum_{n = 1}^\infty ( M(b,a_n,a_{n,l}) + M(b,a_{n+1},a_{n,l}) +  M(a_n,a_{n,l},a_{n,r}) \nonumber \\
& + M(a_{n+1},a_{n,l},a_{n,r}) + M(c,a_n,a_{n,r})+M(c,a_{n+1},a_{n,r}) \nonumber \\
& + M(c,b_n,b_{n,l})+M(c,b_{n+1},b_{n,l})+M(b_n,b_{n,r},b_{n,l}) \nonumber \\
& +M(b_{n+1},b_{n,r},b_{n,l})+ M(a,b_n,b_{n,r})+M(a,b_{n+1},b_{n,r}) \nonumber \\
& + M(a,c_n,c_{n,l})+M(a,c_{n+1},c_{n,l})+ M(c_n,c_{n,l},c_{n,r}) \nonumber \\
& +M(c_{n+1},c_{n,l},c_{n,r})+ M(b,c_n,c_{n,r})+M(b,c_{n+1},c_{n,r}))
\end{align}
(see Figure \ref{mallows necklace}). Since $ t>\delta $,  $ M(a,b,c+\epsilon;t) < M(a,b,c;t)$. Since $ M $ is symmetric in each of $ a,b,c $, it follows that $ M $ is strictly decreasing in each variable.  
\end{proof}

If $ T(a,b,c) $ is dilated by a factor of $\frac{1}{\alpha} $, where $ \alpha > 0  $, then the radii $ r_n $ are replaced by $ \frac{1}{\alpha} r_n $ and $ a,b,c $ are scaled by $ \alpha $. So, $ M $ is homogenous of degree $ -t $:
\begin{equation}\label{M function homogeneity}
M(\al a, \al b, \al c ;t) =  \al^{-t} M(a,b,c;t). 
\end{equation} 

\begin{lemma}\label{basic inequality boyd/malows lemma} Let $ 0 \leq a \leq b \leq c $, $ b>0 $ and $ M $ defined by (\ref{M function B/M packing}). Then, 
	\begin{equation}\label{basic ineq statement b/m}
	(a+c)^{-t}M(0,1,1;t) \leq M(a,b,c;t) \leq b^{-t}M(0,1,1;t)
	\end{equation}  
	with the inequality holding in the extended sense when $ t < S = \dd $. 
\end{lemma}
\begin{remark}
    The Boyd/Mallows packing in $ T(0,1,1) $ is not an integer packing since, for instance $ c_n = 1+2n^2+\sqrt{8}n $. The packing in Figure \ref{boyd mallows packing bounded} contains triangles $T(0,1,2)$ but not $T(0,1,1)$.
\end{remark}

Lemma \ref{basic inequality boyd/malows lemma} may be proved by adopting the methods from Lemma 1 in \cite{Boyd1971disk}. Lemma \ref{basic inequality boyd/malows lemma} shows that $M(a,b,c;t) < \infty $ if and only if
$  M(0,1,1;t) < \infty$, allowing us to analyze the Boyd/Mallows packing within $ T(0,1,1) $. (\ref{basic ineq statement b/m}) shows that $\delta$ is independent of $a,b,c$. Note that since we are assuming that $ 0 \leq a \leq b \leq c $ and $ b>0 $, then $ c_1 \leq b_1 \leq a_1 $, $ a_{n+1} \leq a_{n,l} \leq a_{n,r} $, $ b_{n+1} \leq b_{n,r} \leq b_{n,l} $, and $ c_{n+1} \leq c_{n,l} \leq c_{n,r} $ so that the above sums with $ M(\cdot,\cdot,\cdot) $ in (\ref{functional type equation 0th iteration sum b/m}) are written in increasing curvature. As a result of Lemmas \ref{melzak type lemma}, \ref{left necklace}, and (\ref{functional type equation 0th iteration sum b/m}), a curvature occurring in $T(a,b,c)$ may be written as
\begin{align*}
    \mathrm{curv} = w_1 a + w_2 b + w_3 c + w_4 d
\end{align*}
with $w_j > 0$. Since $M$ is decreasing in all variables, homogenous, and symmetric in $a,b,c$, we get, via \cite{Boyd1973improved}, the following lemma.

\begin{lemma}
    Let $M$ be defined as above and $0 \leq a \leq b \leq c$ and $b>0$. Then,
    \begin{align}\label{improved ineq M function B/M}
        (a+c)^{-t}M(0,1,1;t) \leq \left(a+\frac{b+c}{2} \right)^{-t} M(0,1,1;t) \leq M(a,b,c;t), \nonumber \\ M(a,b,c;t) \leq \frac{1}{2}((a+b)^{-t}+c^{-t})M(0,1,1;t) \leq b^{-t}M(0,1,1;t).
    \end{align}
\end{lemma}

To help us understand the forthcoming definitions, we will first exhibit how to create a self-similar inequality for $M$. Based on (\ref{functional type equation 0th iteration sum b/m}), define
\begin{align*}
   h_0(t) =  \sum_{n=1}^\infty \left( a_n^{-t} + a_{n,l}^{-t} + a_{n,r}^{-t} + b_n^{-t} + b_{n,l}^{-t} + b_{n,r}^{-t} + c_n^{-t} + c_{n,l}^{-t} + c_{n,r}^{-t}  \right). 
\end{align*}
Applying the rightmost inequality of Lemma \ref{basic inequality boyd/malows lemma}, we get
\begin{align}
& M(a,b,c) \leq h_0(t) + 2M(0,1,1)(c_1^{-t}+b_1^{-t} \nonumber \\
& + \sum_{n=1}^\infty (a_n^{-t}+a_{n+1}^{-t}+a_{n,l}^{-t} + b_n^{-t}+b_{n+1}^{-t}
+b_{n,r}^{-t} + c_n^{-t}+c_{n+1}^{-t}+c_{n,l}^{-t})). 
\end{align}
If we write
\begin{align*}
    \tilde{f}_0(t) & = 2(c_1^{-t}+b_1^{-t} \\
    & + \sum_{n=1}^\infty (a_n^{-t}+a_{n+1}^{-t}+a_{n,l}^{-t} + b_n^{-t}+b_{n+1}^{-t}+b_{n,r}^{-t} + c_n^{-t}+c_{n+1}^{-t}+c_{n,l}^{-t}))
\end{align*}
and let $ (a,b,c) = (0,1,1) $, then
\begin{equation}
M(0,1,1;t) \leq h_0(t) + M(0,1,1;t)\tilde{f}_0(t).
\end{equation}
If $\tilde{\mu}_0$ satisfies $\tilde{f}_0(\tilde{\mu}_0) = 1$, then Theorem \ref{theorem 1 b/m} will show that $\tilde{\mu}_0 > S$. Applying the leftmost inequality from Lemma \ref{basic inequality boyd/malows lemma} to equation (\ref{functional type equation 0th iteration sum b/m}) gives
\begin{align}
M(a,b,c) & \geq h_0(t)  \nonumber \\ & + M(0,1,1)\big((a+b_1)^{-t}+(b+a_1)^{-t}+(c+a_1)^{-t}+(c_1+a_1)^{-t} \nonumber \\
& + \sum_{n=1}^\infty (2(b+a_{n,l})^{-t}+ 2(c+a_{n,r})^{-t} + (a_n+a_{n,r})^{-t} + (a_{n+1}+a_{n,r})^{-t}  \nonumber \\
& + 2(c+b_{n,l})^{-t} + 2(a+b_{n,r})^{-t} + (b_n+b_{n,l})^{-t} + (b_{n+1}+b_{n,l})^{-t} \nonumber \\
& + 2(a+c_{n,l})^{-t}+2(b+c_{n,r})^{-t}+(c_n+c_{n,r})^{-t} + (c_{n+1}+c_{n,r})^{-t}) \big). 
\end{align}
Letting $ \tilde{g}_0(t) $ be the function to the right of $ M(0,1,1) $ above with $ (a,b,c) = (0,1,1) $, we can write
\[ M(0,1,1) \geq h_0(t) + M(0,1,1)\tilde{g}_0(t). \]
Theorem \ref{theorem 1 b/m} will also show that if $ \tilde{g}_0(\tilde{\lambda}_0) = 1 $, then $\tilde{\lambda}_0 \leq S$. 

We now define functions similar to $\tilde{f}_0,\tilde{g}_0$, and $h_0$, but based on the tighter inside inequalities of (\ref{improved ineq M function B/M}). First define a set-valued function $\tau$ which collects the triples of curvatures $(x,y,z)$ when $T(x,y,z)$ occurs as a sub-triangle per Figure \ref{mallows necklace}. Define
\begin{align}\label{triangle_list function}
\tau(a,b,c) &  = \{ (a,c_1,b_1),(b,c_1,a_1),(c,b_1,a_1),(c_1,b_1,a_1)  \} \cup \nonumber \\
& \bigcup_{n=1}^\infty \{ (b,a_n,a_{n,l}) , (b,a_{n+1},a_{n,l}) ,  (a_n,a_{n,l},a_{n,r}),  \nonumber \\
& (a_{n+1},a_{n,l},a_{n,r}), (c,a_n,a_{n,r}), (c,a_{n+1},a_{n,r}),  \nonumber \\
& (c,b_n,b_{n,l}),(c,b_{n+1},b_{n,l}),(b_n,b_{n,r},b_{n,l}), \nonumber \\
& (b_{n+1},b_{n,r},b_{n,l}),(a,b_n,b_{n,r}),(a,b_{n+1},b_{n,r}), \nonumber \\
& (a,c_n,c_{n,l}),(a,c_{n+1},c_{n,l}),(c_n,c_{n,l},c_{n,r}), \nonumber \\
& (c_{n+1},c_{n,l},c_{n,r}),(b,c_n,c_{n,r}),(b,c_{n+1},c_{n,r})  \}. 
\end{align}
Define
\begin{align}
& \tilde{f}_0(\kk;a,b,c;t) = \sum_{(x,y,z) \in{\tau(a,b,c)}}  y^{-t},  \\
& \tilde{g}_0(\kk;a,b,c;t) = \sum_{(x,y,z)\in{\tau(a,b,c)}}  (x+z)^{-t}, \\
& f_0(\kk;a,b,c;t) = \sum_{(x,y,z) \in{\tau(a,b,c)}}  \frac{1}{2}((x+y)^{-t}+z^{-t}),  \\
& g_0(\kk;a,b,c;t) = \sum_{(x,y,z)\in{\tau(a,b,c)}}  \left(x+\frac{y+z}{2}\right)^{-t}, \\
& h_0(\kk;a,b,c;t) =  \sum_{n=1}^\infty \left( a_n^{-t} + a_{n,l}^{-t} + a_{n,r}^{-t} + b_n^{-t} + b_{n,l}^{-t} + b_{n,r}^{-t} + c_n^{-t} + c_{n,l}^{-t} + c_{n,r}^{-t}  \right). 
\end{align}
For $\kappa \geq 0$, define 
\begin{align}\label{list_iter function}
\mathscr{S}(\kk;\tau(a,b,c)) & = \tau(a,b,c) \cup \bigcup _{\substack{(x,y,z)\in{\tau(a,b,c)}\\
                  y < \kappa}} \tau(x,y,z) \nonumber \\
            & - \lbrace (p,q,r)\in{\tau(a,b,c)}: q < \kappa \rbrace
\end{align}
and $\mathscr{S}^0( \kappa;l) = l$. Define iterates of $\mathscr{S}$ as
\begin{align*}
    \mathscr{S}^2(\kappa;l) & = \mathscr{S}(\kappa; \mathscr{S}(\kappa;l)) \\
    \mathscr{S}^3(\kappa;l) & = \mathscr{S}(\kappa; \mathscr{S}(\kappa;\mathscr{S}(\kappa;l)))
\end{align*}
etc. For $ m \geq 1 $, define
\begin{align} \label{fm function B/M}
\tilde{f}_m(\kk;a,b,c;t) &  = \sum_{(x,y,z)\in{\mathscr{S}^m(\kappa;\tau(a,b,c))}} y^{-t},  \\
\tilde{g}_m(\kk;a,b,c;t) &  = \sum_{(x,y,z)\in{\mathscr{S}^m(\kappa;\tau(a,b,c))}} (x+z)^{-t},  \\
f_m(\kk;a,b,c;t) &  = \sum_{(x,y,z)\in{\mathscr{S}^m(\kappa;\tau(a,b,c))}} \frac{1}{2}((x+y)^{-t}+z^{-t})  \\
g_m(\kk;a,b,c;t) & = \sum_{(x,y,z)\in{\mathscr{S}^m(\kappa;\tau(a,b,c))}} \left(x+\frac{y+z}{2}\right)^{-t} \label{gm function B/M}   \\
h_m(\kk;a,b,c;t) & = h_{m-1}(\kk;a,b,c;t) \nonumber \\
& + \sum_{\substack{(x,y,z)\in{\mathscr{S}^{m-1}(\kappa;\tau(a,b,c))}\\ y < \kappa}} h_0(\kk;x,y,z;t) \label{hm function B/M} 
\end{align}
For fixed $ \kk,a,b,c$, and $m$, the functions $ f_m,g_m, $ and $ h_m $ above are defined when $ t > \frac{1}{2} $. At $ t = \frac{1}{2} $, they are harmonic-type series since $ c_n \asymp n^2 $. If $ c_1 > 1 $, they are positive, continuous, strictly monotone decreasing functions of $ t $, tending to $ \infty $ as $ t \rightarrow \frac{1}{2}^+ $ and $ 0 $ as $ t \rightarrow \infty $. Based on the above definitions of $ f_m,g_m, $ and $ h_m $, if $ \kk $ is smaller than all sub-triangle curvatures, we expect that the iteration should terminate; in other words, $ f_m=f_{m+1}= \cdots  $.  

\begin{figure}[h!]
	\begin{center}
		\includegraphics[scale=0.4]{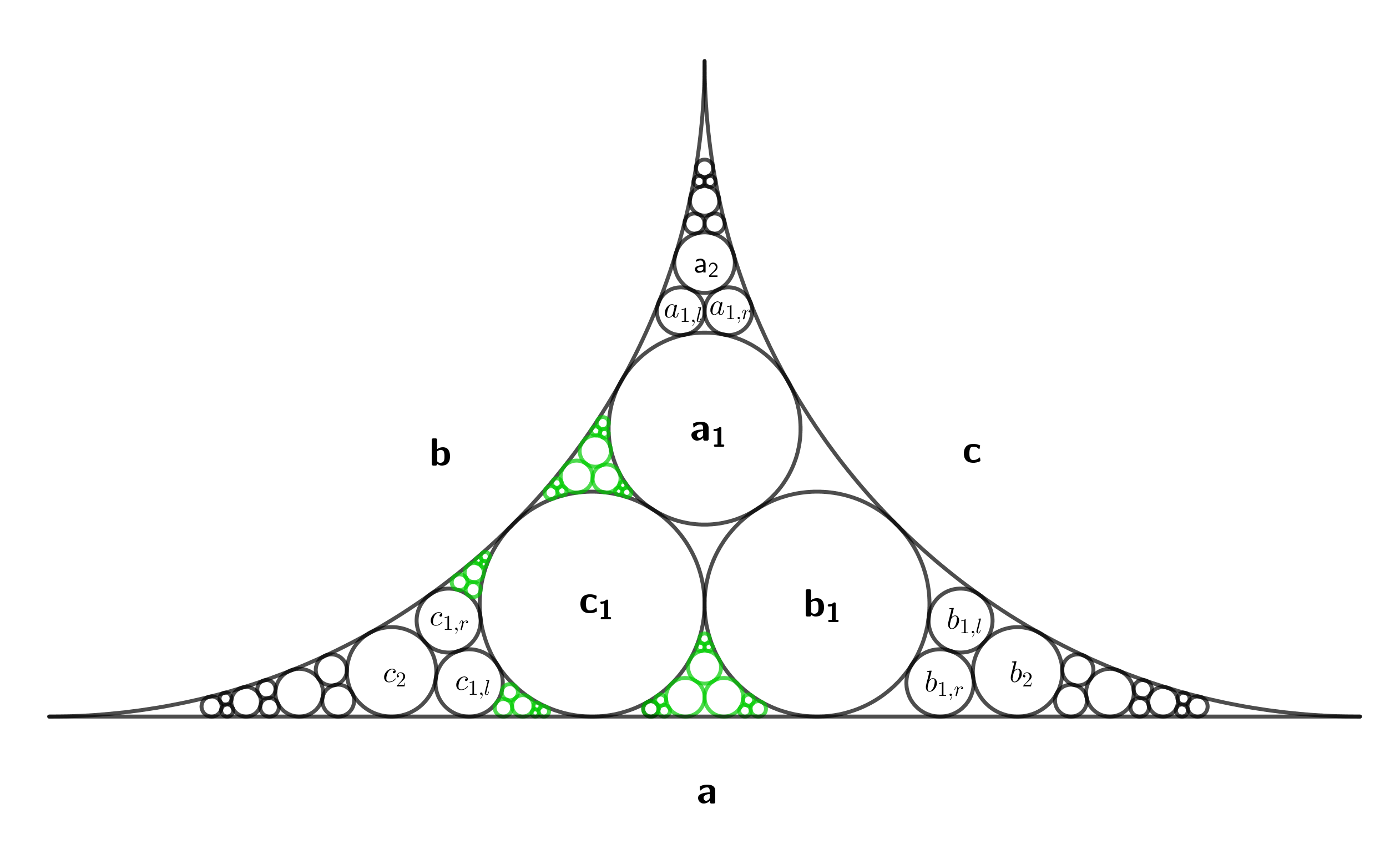}
	\end{center}
	\caption{A first iteration of the necklace packing for any $ c_1 < \kk \leq b_1 $. The green/gray disks 
		are added because those four sub-triangles have a middle curvature smaller than $ \kk $.}
\end{figure}
The following lemma establishes that the breaking of necklace triangles into sub-triangles occurs with curvatures which grow at an exponential rate.

\begin{lemma}
	If $ \kk \leq 5^{n+1}b $, then
	\begin{align*}
	\mathscr{S}^{n+1}(\kappa,\tau(a,b,c)) = \mathscr{S}^n(\kappa,\tau(a,b,c)) \quad \text{ for all } n \geq 0.
	\end{align*}
\end{lemma}
    \begin{proof}:
        On any triangle triple $(a,b,c)$ with $0 \leq a \leq b \leq c$ and $b>0$, we have 
    \begin{align*}
        \min \{ y: (x,y,z)\in{\tau(a,b,c)} \} = c_1
    \end{align*}
    and
    \begin{align*}
       c_1 = c+2a+2b+\sqrt{8}\sqrt{ab+ac+bc} \geq 3b + \sqrt{8}\sqrt{b^2} > 5b. 
    \end{align*}
    So, if $\kappa \leq 5b < c_1$, then by (\ref{list_iter function}),
    \begin{align*}
    \mathscr{S}(\kappa;\tau(a,b,c)) = \mathscr{S}^0(\kappa;\tau(a,b,c)) = \tau(a,b,c).    
    \end{align*}
    If true for $1,...,n-1$, then by choosing $\kappa \leq 5^{n+1}b$, the union and set exclusion on the right hand side of (\ref{list_iter function}) for $\mathscr{S}^{n+1}$ will be empty, and the claim follows.
    \end{proof}

\begin{corollary} \label{boyd lemma 4}
	Fix $n \geq 0$ and let $ m \geq 0 $. If $ \kk \leq 5^{n+1}b $ then
	\begin{align*}
	f_m(\kk;a,b,c;t) &= f_{n}(\kk;a,b,c;t) \quad (m \geq n) \\
	g_m(\kk;a,b,c;t) &= g_{n}(\kk;a,b,c;t) \quad  (m \geq n) \\
	h_m(\kk;a,b,c;t) &= h_{n}(\kk;a,b,c;t) \quad  (m \geq n).
	\end{align*}
\end{corollary}
When $ (a,b,c) = (0,1,1) $, the base case $ n=0 $ gives $ c_1 = 3+2\sqrt{2} $ and by the symmetry of the packing in $T(0,1,1)$,
	\[  \min \{y:(x,y,z) \in{\mathscr{S}^m(\infty;\tau(0,1,1))} \} = (3+2\sqrt{2})^{m+1}.   \]
The $\kappa$-cutoff values in Table \ref{table 1 upper and lower bounds b/m} are $ \lfloor \beta_m \rfloor$ where $\beta_m = \beta_0^{m+1}$, $ \beta_0 =  3+2 \sqrt{2}$.
\begin{theorem}\label{theorem 1 b/m}
	Let $ S$ be the critical exponent of the Boyd/Mallows packing. Let $ \kk > 0 $, $ m \geq 0 $, and $ t > \frac{1}{2} $. Let $ \mu_m(\kk) $ and $ \lambda_m(\kk) $ satisfy $ f_m(\kk;0,1,1;\mu_m(\kk))  = 1$ and $ g_m(\kk;0,1,1;\lambda_m(\kk)) = 1 $. Then,
	\begin{equation}
	\lambda_m(\kk) \leq S \leq \mu_m(\kk).
	\end{equation}
	Moreover, if $ 1 < \kk \leq (3+2 \sqrt{2})^{m+1} $, then   
	\begin{equation}\label{boyd 27}
	0 < \mu_m(\kappa) - \lambda_m(\kk) < \frac{2.3}{\log \kk}.
	\end{equation}
\end{theorem}

\begin{proof}:
	With $ (a,b,c) = (0,1,1) $, $ c_1 = 3+2\sqrt{2} > 1 $, so $ f_m$ and $g_m $ are  strictly decreasing, continuous functions. Thus, $ \mu_m(\kk) $ and $ \lambda_m(\kk) $ are unique. Define
    \begin{align}\label{finite M function}
    M_j(a,b,c;t) = \sum_{q \leq c_1^j} q^{-t},
    \end{align}
    where the sum is over curvatures $q$ occurring in $T(a,b,c)$. We will first show that
	\begin{equation}\label{boyd 30}
	M_j(0,1,1;t) \leq h_m(\kk;0,1,1;t) + M_j(0,1,1;t)f_m(\kk;0,1,1;t). 
	\end{equation}	
	Let $ m = 0 $, fix $ \kk > 0 $, and $ t > \frac{1}{2} $. Using (\ref{improved ineq M function B/M}) applied to $M_j$, 
	\begin{align}
	M_j(a,b,c;t) & < h_0(\kk;a,b,c;t) + \sum_{\substack{(x,y,z)\in{\tau(a,b,c)}\\
                  }} M_j(x,y,z;t)  \nonumber \\
    & \leq h_0(\kappa;a,b,c;t) + \sum_{(x,y,z)\in{\tau(a,b,c)}} \frac{1}{2}((x+y)^{-t}+z^{-t})M_j(0,1,1;t) \nonumber \\
	& = h_0(\kk;a,b,c;t) + f_0(\kappa;a,b,c;t)M_j(0,1,1;t). \nonumber \\
	\end{align}
	If also true for $ 1,...,m-1 $, then using the definitions of $ f_m $ and $ h_m $, 
	\begin{align}
	M_j(a,b,c;t) \quad \quad &  \nonumber \\
	< h_0(\kk;a,b,c;t) & + \sum_{\substack{(x,y,z)\in{\tau(a,b,c)}\\ y < \kappa}} M_j(x,y,z;t)  \nonumber \\
	& + \sum_{\substack{(x,y,z)\in{\tau(a,b,c)}\\ \kappa \leq y }} M_j(x,y,z;t) \nonumber \\  
	< h_0(\kk;a,b,c;t) & \nonumber \\
    + \sum_{\substack{(x,y,z)\in{\tau(a,b,c)}\\ y < \kappa}} &(h_{m-1}(\kappa;x,y,z;t)+f_{m-1}(\kappa;x,y,z;t)M_j(0,1,1;t)) \nonumber \\ 
	 + \sum_{\substack{(x,y,z)\in{\tau(a,b,c)}\\ \kappa \leq y}} & \frac{1}{2}((x+y)^{-t}+z^{-t}) M_j(0,1,1;t) \nonumber \\
	< h_m(\kk;a,b,c;t&)  + M_j(0,1,1;t)f_m(\kk;a,b,c;t).
	\end{align}
	This establishes (\ref{boyd 30}) by induction. Now let $ t>\mu_m(\kk) $, making $ f_m(\kk;0,1,1;t) < 1 $. Then,
	\begin{align}
	M_j(0,1,1;t) < \frac{h_m(\kk;0,1,1;t)}{1-f_m(\kk;0,1,1;t)}.
	\end{align}
	Since $ M_j \nearrow M $, letting $ j \rightarrow \infty $, shows that when $ t > \mu_m(\kk) $, $ M(0,1,1;t) $ is bounded above by the finite quantity $h_m(\kappa;0,1,1;t)/(1-f_m(\kappa;0,1,1;t)$. Thus $ \mu_m(\kk) \geq S $. Using (\ref{functional type equation 0th iteration sum b/m}) and (\ref{improved ineq M function B/M}), we may write $ M(a,b,c;t) $ as follows:
	\begin{align*}
	M(a,b,c;t) & = h_m(\kk;a,b,c;t) + \sum_{(x,y,z)\in{\mathscr{S}^m(\kappa; \tau(a,b,c))}} M(x,y,z;t) \\
	& \geq h_m(\kk;a,b,c;t) + \sum_{(x,y,z)\in{\mathscr{S}^m(\kappa; \tau(a,b,c))}} \left(a+\frac{b+c}{2}\right)^{-t} M(0,1,1;t) \\
	& = h_m(\kk;a,b,c;t) + g_m(\kk;a,b,c;t)M(0,1,1;t).
	\end{align*}
	Letting $ (a,b,c) = (0,1,1) $ establishes
 \begin{align*}
     \frac{h_m(\kappa;0,1,1;t)}{1-g_m(\kappa;0,1,1;t)} \leq M(0,1,1;t).
 \end{align*}
    If $ t > \lambda_m(\kk) $ then $ g_m < 1 $ since $ g_m $ is decreasing. Since 
    \begin{align*}
      \lim_{t \rightarrow \lambda_m(\kk)^+}g_m(\kk;0,1,1;t) = 1  
    \end{align*}
    and  
	\begin{align*}
	\lim_{t \rightarrow \lambda_m(\kk)^+} \frac{h_m(\kk;0,1,1;t)}{1-g_m(\kk;0,1,1;t)} \leq \lim_{t \rightarrow \lambda_m(\kk)^+} M(0,1,1;t) = M(0,1,1;\lambda_m(\kk))
	\end{align*}
	then $ M(0,1,1;\lambda_m(\kk)) = \infty $, which shows that $ \lambda_m(\kk) \leq S $.
	Next, we show that for fixed $\kappa,a,b,c,t$, 
	\begin{align}\label{boyd 36}
	g_m \geq  \tilde{g}_m \geq 5.5^{-t}\tilde{f}_m \geq 5.5^{-t}f_m 
	\end{align}
	The outermost inequalities follow from (\ref{improved ineq M function B/M}). For the inside inequality, with $ m = 0 $, we compare $ x+z $ to $ y $ for $ (x,y,z) \in{\tau(a,b,c)} $. For instance, if $ (x,y,z) = (c,a_n,a_{n,r}) $ (the $ 9^\text{th} $ necklace term in $ \tau(a,b,c) $) then
	\begin{align*}
	c+a_{n,r} & = 2a+b(1+4n+4n^2)+c(3+4n+4n^2)+\sqrt{8}d(2n+1) \\
	& \leq 5.5( a+b(2n^2) + c(2n^2) + \sqrt{8}dn) = 5.5a_n.  
	\end{align*}
	The $ 11^{\text{th}} $ and $ 21^{\text{st}} $ terms in $ \tau(a,b,c) $ can also be compared with the constant $ 5.5 $, while the other terms can be compared with constants (also independent of $ \kk $ and $ m $) ranging between $ 2 $ and $ 5 $. So, $ (x+z)^{-t} \geq 5.5^{-t}y^{-t} $ for all $ (x,y,z) \in{\tau(a,b,c)} $ implying $ \tilde{g}_0(\kk;a,b,c;t) \geq 5.5^{-t}\tilde{f}_0(\kk;a,b,c;t) $. By the definition of $\mathscr{S}$ and induction on $m$, we get that $(x+z)^{-t} \geq 5.5^{-t}y^{-t}$ for all $(x,y,z)\in{\mathscr{S}^m(\kappa;\tau(a,b,c))}$, establishing (\ref{boyd 36}).
	
	Now let $ 1 < \kk \leq \beta_m $. By Corollary \ref{boyd lemma 4}, $ \lambda_m(\kk) = \lambda_{m+1}(\kk) = ... $ and $ \mu_m(\kk) = \mu_{m+1}(\kk) = ... $ and $ y \geq \kk $ for all $ (x,y,z) \in{\mathscr{S}^m(\kk;\tau(a,b,c))} $. Let $ \epsilon > 0 $ be given. Since $ \frac{\kk}{\beta_m} \leq 1 $, then if $ q $ is a curvature occurring in the expression of $ \tilde{f}_m(\kk;0,1,1;t) $, then $ q \geq \beta_m $, so $ \frac{\kk}{q} \leq 1 $. Thus, $ \left( \frac{\kk}{q}  \right)^\epsilon \leq 1^\epsilon = 1 $ and so 
	\begin{align}\label{multiplicative inequality f b/m}
	\tilde{f}_m(\kk;0,1,1;t) = \sum_{n=1}^\infty q_n^{-t} \geq \sum_{n = 1}^\infty \kk^\epsilon q_n^{-t-\epsilon} = \kk^\epsilon \tilde{f}_m(\kk;0,1,1;t+\epsilon).
	\end{align}
	Numerical computation (see next section) showed that $ \mu_5(39201) = 1.348771 $, so $ \lambda_m(\kk) \leq 1.348771 $ for all $ \kk $ and $ m $. Thus, $ 5.5^{-\lambda_m(\kk)} \geq 5.5^{-1.348771} = 0.100327$. Using (\ref{boyd 36}), (\ref{multiplicative inequality f b/m}), setting $\mathcal{C} = 0.100327$, $\epsilon = \mu_m(\kappa)-\lambda_m(\kappa)$, and $ \kk>1 $ to ensure $ \log \kk > 0  $,  
	\begin{align} \label{convergence upper lower b/m}
	1 & = g_m(\kk;0,1,1;\lambda_m{(\kk)}) \nonumber \\
	& \geq 5.5^{-\lambda_m(\kk)}\tilde{f}_m(\kk;0,1,1;\lambda_m(\kk)) \nonumber \\
	& \geq \mathcal{C} \kappa^\epsilon f_m(\kk;0,1,1;\lambda_m(\kk)+\epsilon) \nonumber \\
	& = \mathcal{C} \kk^{\mu_m(\kk) - \lambda_m(\kk)} f_m(\kk;0,1,1; \mu_m(\kk)) \nonumber \\
	& = \mathcal{C} \kk^{\mu_m(\kk) - \lambda_m(\kk)}. 
	\end{align} 
	This proves (\ref{boyd 27}) and concludes the proof of the theorem. 
\end{proof}

\begin{remark}
The above condition $ 1 < \kk \leq \beta_m $ is necessary in establishing the multiplicative inequality (\ref{multiplicative inequality f b/m}). For example, if $ m=0 $ and $ \kk = 6 > \beta_0 = 3 + \sqrt{8} $, then $ \tau(0,1,1) $ includes $ (0,3+\sqrt{8},3+\sqrt{8}) $ so the property $ \kk / q \leq 1 $ fails. If $ 0 < \kk \leq 1 $, then $ \log \kk \leq 0 $ and the inequalities leading to (\ref{convergence upper lower b/m}) flip and give no useful information. If $ \kk = 1 $, (\ref{multiplicative inequality f b/m}) degrades to the (already known) monotonic decreasing property of $ \tilde{f}_m $. $ 0 < \kk \leq 1 $ could be treated as a separate case, but since in our computations we seek values of $ \kk > 10^4 $, the insistence that $ \kk > 1 $ is a practical condition. 
\end{remark}

\section{Numerical Results}

The functions $\tau, \mathscr{S}, f_m$, and $g_m$ were incorporated into a computer program written using Sage mathematics software \cite{stein2007sage}. All computations were run on personal computers with $ 8 $ Gigabytes of random access memory and a single central processor unit. Table \ref{table 1 upper and lower bounds b/m} provides the computed upper and lower bounds. The convergence rate with both sets of inequalities appears to be roughly $ \frac{1}{\log \kk} $. The functions $ \lambda_{m} $ and $ \mu_{m} $ (which are based on (\ref{improved ineq M function B/M})) provide closer starting values ($ \lambda_{0}(\kk) $ and $ \mu_{0}(\kk) $) to $ S $ as compared to $\tilde{\lambda}_m$ and $\tilde{\mu}_m$ (which are based on (\ref{basic ineq statement b/m})). As suggested by Theorem \ref{theorem 1 b/m}, the approximate slope of $ f_m,g_m,\tilde{f}_m,$ and $\tilde{g}_{m} $ is $ -\log \kappa $. A numerical root-finding algorithm incorporating a Newton type map $ x_0 \mapsto -\frac{A}{\log \kk}(1-f_0(\kk;0,1,1,x_0)) + x_0 $ was employed using a constant approximate slope. The value of $ A \approx 1.5 $ as well as the initial guess for $ x_0 $ can be adjusted to offer faster or slower convergence to the root of $ f_m-1 $, $ g_m-1 $, $ \tilde{f}_{m}-1 $, or $ \tilde{g}_{m}-1 $. Iterations of the Newton type map were done until $ |f_m(\kk;0,1,1;x_n) - 1| < 10^{-7} $,  $ x_n = x_0,x_1,x_2,... $ and similar for $ g_m $, $ \tilde{f}_{m} $, $ \tilde{g}_{m} $.  

\begin{table}[h!]
	\centering
	\begin{tabular}{|c c c c c c|} 
        \hline
		$ m $ & $ \kk $ & $ \tilde{\lambda}_m(\kk) $ & $ \lambda_{m}(\kk) $ & $ \mu_{m}(\kk) $ & $ \tilde{\mu}_m(\kk) $ 
		\\ [0.5ex] 
		\hline \hline
		0 &  $\leq \beta_0$  & 1.238656 & 1.304679 &  1.391406  & 1.549702\\ 
		1 & 16 & 1.274746 & 1.316674 & 1.367061 & 1.445461 \\
		1 & 33 & 1.278722 & 1.318153 & 1.365074 & 1.437800 \\
		2 & 100 & 1.288116 & 1.320996 & 1.359760 & 1.417712 \\
		2 & 197 & 1.292704 & 1.322415 & 1.357262 & 1.408436 \\
		3 & 1153 & 1.300423 & 1.324607 & 1.353417 & 1.394080 \\
		4 & 6725 & - & 1.326166 & 1.350711 & - \\
		5 & 39201 & - & 1.327266 & 1.348771 & - \\
		\hline
	\end{tabular}
	\caption{Table of upper and lower bounds of $ S = \dd $ for various $ \kk $ and $ m $ of the Boyd/Mallows packing.}
	\label{table 1 upper and lower bounds b/m}
\end{table}

The $ \kk $ values of $ 33, 197,  1153, 6725 $, and $ 39201 $ were rounded down from the cutoff values of $ \beta_1 = \beta_0^2 \approx 33.97$, $ \beta_2 = \beta_0^3 \approx 197.99 $ etc.,  where $ \beta_0 = 3 + 2 \sqrt{2} \approx 5.82 $. The other values of $ 16$ and $ 100 $ were inserted to illustrate that further refinement can happen even with the same number of iterations. The fact that $ \lambda_{1}(16) = 1.316674 > 1.310876 > \delta_{\mathcal{A}} $ proves that the critical exponent of the Boyd/Mallows packing is strictly greater than the critical exponent of the Apollonian packing. The longest calculations took approximately 6 hours to complete. The above values are truncated, not rounded. 

To obtain bounds for the critical exponent of the Apollonian packing, $\delta_{\mathcal{A}}$, analogous functions $f_{m, \mathcal{A}},g_{m,\mathcal{A}},\lambda_{m,\mathcal{A}},\mu_{m,\mathcal{A}},$ were created. A $\kappa$-cutoff value of $166464$ was used to provide 6 levels of iteration. Obtaining the roots $\lambda_{6,\mathcal{A}}$ and $\mu_{6,\mathcal{A}}$ yielded (\ref{apollonian new bound}).

We next obtain a heuristic estimate of $ S = \dd $ by methods which originate in \cite{Boyd1973Separation}, \cite{Melzak1966Infinite}, and \cite{Gilbert1964Randomly}. Since each generator of $ \mathcal{O}_J^+(\ZZ)$ has a matrix representation, the problem of estimating $\delta$ is equivalent to determining the average growth rate of the orbit of a vector under random products of certain non-commuting matrices. After a similarity transformation, these integer matrices  correspond to Lorentz transformations (see \cite{Thomas1994Hausdorf} or \cite{Soderberg1992apollonian}) which generate a discrete subgroup of the full Lorentz group $ O_{1,3}(\RR) $. Define the ``height" function
\begin{align}\label{b/m height function}
h(\vx) = \vx \circ \vD = 4x_1+4x_2+2x_3+2x_4.
\end{align}
Again let $\Gamma = \la R_{\vn_1},R_{\vn_2},R_{\vn_3},R_{\vn_4},R_{\vn_5}  \ra$ and consider the image of $\Gamma$ on the faces $\{ \ve_1,\ve_2,\ve_3,\ve_4,\ve_1+\ve_2-\ve_4 \}$ (see Figure \ref{b/m symmetries}) . A finite subset of this image was created using the Sagemath function
\texttt{RecursivelyEnumeratedSet} with the 5 faces \texttt{[e1,e2,e3,e4,e1+e2-e4]} as the \texttt{seed} and any \texttt{successor} defined to be any vector in the orbit with a maximum height \texttt{hmax} below $ 2^{19} $. For improved memory management, the option \texttt{structure = `symmetric'} was used since each generator is a reflection. There were calculated to be $ 13,244,370 $ vectors with a height below $ 2^{19} $. These vectors were then sorted according to height and then fit by linear least-squares regression to the curve $ y = ax^b $ using the function \texttt{find fit}. The resulting heuristic estimate of $ 1.33544546879 $ fits between the rigorous bounds displayed in Table \ref{table 1 upper and lower bounds b/m} and is merely $ 0.002573 $ off from the arithmetic average of the lower and upper bounds $ \lambda_{5}(39201) $ and $ \mu_{5}(39201) $.

\section*{Acknowledgments}
I would like to thank Arthur Baragar for many helpful conversations. Also thanks to David Lautzenheiser and Judson Clark for help in improving my computer code.
\newpage

\newpage
\bibliographystyle{plain}

\bibliography{references}

\begin{thebibliography}{10}

\bibitem{Apanasov2011Conformal}
Boris~N Apanasov.
\newblock {\em Conformal geometry of discrete groups and manifolds}, volume~32.
\newblock Walter de Gruyter, 2011.

\bibitem{Baragar2001Survey}
Arthur Baragar.
\newblock A survey of classical and modern geometries.
\newblock {\em Upper Saddle River}, 2001.

\bibitem{Baragar2017Higher}
Arthur Baragar.
\newblock Higher dimensional apollonian packings, revisited.
\newblock {\em Geometriae Dedicata}, pages 1--25, 2017.

\bibitem{baragar2020game}
Arthur Baragar and Daniel Lautzenheiser.
\newblock A game of packings.
\newblock {\em arXiv preprint arXiv:2006.00619}, 2020.

\bibitem{Boyd1974New}
David Boyd.
\newblock A new class of infinite sphere packings.
\newblock {\em Pacific Journal of Mathematics}, 50(2):383--398, 1974.

\bibitem{Boyd1971disk}
David~W Boyd.
\newblock The disk-packing constant.
\newblock {\em aequationes mathematicae}, 7(2):182--193, 1971.

\bibitem{Boyd1973improved}
David~W Boyd.
\newblock Improved bounds for the disk-packing constant.
\newblock {\em Aequationes Mathematicae}, 9(1):99--106, 1973.

\bibitem{Boyd1973Separation}
David~W Boyd.
\newblock The osculatory packing of a three-dimensional sphere.
\newblock {\em Canadian J. Math}, 25:303--322, 1973.

\bibitem{Clifford1882Powers}
WK~Clifford.
\newblock On the powers of spheres (1868).
\newblock {\em Mathematical Papers of William Kingdon Clifford}, pages 332--336, 1882.

\bibitem{Darboux1872Points}
MG~Darboux.
\newblock De points, de cercles et de spheres.
\newblock In {\em Annales de L'Ecole Normale}, volume~1, pages 323--392, 1872.

\bibitem{Dolgachev2016Orbital}
Igor Dolgachev.
\newblock Orbital counting of curves on algebraic surfaces and sphere packings.
\newblock In {\em K3 Surfaces and Their Moduli}, pages 17--53. Springer, 2016.

\bibitem{Gilbert1964Randomly}
EN~Gilbert.
\newblock Randomly packed and solidly packed spheres.
\newblock {\em Canad. J. Math}, 16:286--298, 1964.

\bibitem{Mallows2010Generalization}
Gerhard Guettler and Colin Mallows.
\newblock A generalization of apollonian packing of circles.
\newblock {\em Journal of Combinatorics}, 1(1):1--27, 2010.

\bibitem{kontorovich2011apollonian}
Alex Kontorovich and Hee Oh.
\newblock Apollonian circle packings and closed horospheres on hyperbolic 3-manifolds.
\newblock {\em Journal of the American Mathematical Society}, 24(3):603--648, 2011.

\bibitem{Manna1991Precise}
SS~Manna and HJ~Herrmann.
\newblock Precise determination of the fractal dimensions of apollonian packing and space-filling bearings.
\newblock {\em Journal of Physics A: Mathematical and General}, 24(9):L481, 1991.

\bibitem{Maxwell1981Space}
George Maxwell.
\newblock Space groups of coxeter type.
\newblock In {\em Proceedings of the Conference on Kristallographische Gruppen (Univ. Bielefeld, Bielefeld, 1979), Part II}, pages 65--76, 1981.

\bibitem{Melzak1966Infinite}
ZA~Melzak.
\newblock Infinite packings of disks.
\newblock {\em Canad. J. Math}, 18:838--853, 1966.

\bibitem{Mergelyan1952uniform}
Sergei~Nikitich Mergelyan.
\newblock Uniform approximations of functions of a complex variable.
\newblock {\em Uspekhi Matematicheskikh Nauk}, 7(2):31--122, 1952.

\bibitem{Ratcliffe2006Foundations}
John Ratcliffe.
\newblock {\em Foundations of hyperbolic manifolds}, volume 149.
\newblock Springer Science \& Business Media, 2006.

\bibitem{Soderberg1992apollonian}
Bo~S{\"o}derberg.
\newblock Apollonian tiling, the lorentz group, and regular trees.
\newblock {\em Physical Review A}, 46(4):1859, 1992.

\bibitem{stein2007sage}
William Stein.
\newblock Sage mathematics software.
\newblock {\em http://www. sagemath. org/}, 2007.

\bibitem{Sullivan1984Entropy}
Dennis Sullivan.
\newblock Entropy, hausdorff measures old and new, and limit sets of geometrically finite kleinian groups.
\newblock {\em Acta Mathematica}, 153(1):259--277, 1984.

\bibitem{Thomas1994Hausdorf}
PB~Thomas and Deepak Dhar.
\newblock The hausdorf dimension of the apollonian packing of circles.
\newblock {\em Journal of Physics A: Mathematical and General}, 27(7):2257, 1994.

\bibitem{wesler1960infinite}
Oscar Wesler.
\newblock An infinite packing theorem for spheres.
\newblock {\em Proceedings of the American Mathematical Society}, 11(2):324--326, 1960.

\end{thebibliography}

\end{document}